\documentstyle[12pt]{article}

\renewcommand{\thefootnote}{\fnsymbol{footnote}}
\makeatletter \ifcase\@ptsize \font\teneufm=eufm10
\font\seveneufm=eufm7 \font\fiveeufm=eufm5 \font\teneusm=eusm10
\font\seveneusm=eusm7 \font\fiveeusm=eusm5 \or
\font\teneufm=eufm10 scaled \magstephalf \font\seveneufm=eufm7
\font\fiveeufm=eufm5 \font\teneusm=eusm10 scaled \magstephalf
\font\seveneusm=eusm7 \font\fiveeusm=eusm5 \or
\font\teneufm=eufm10 scaled \magstep1 \font\seveneufm=eufm7
\font\fiveeufm=eufm5 \font\teneusm=eusm10 scaled \magstep1
\font\seveneusm=eusm7 \font\fiveeusm=eusm5 \fi

\newfam\eufmfam \newfam\eusmfam \textfont\eufmfam=\teneufm
\scriptfont\eufmfam=\seveneufm \scriptscriptfont\eufmfam=\fiveeufm
\textfont\eusmfam=\teneusm \scriptfont\eusmfam=\seveneusm
\scriptscriptfont\eusmfam=\fiveeusm

\def\frak{\ifmmode\let\next\frak@\else
 \def\next{\errmessage{Use \string\frak\space only in math
 mode}}\fi\next} \def\frak@#1{{\frak@@{#1}}}
 \def\frak@@#1{\fam\eufmfam#1} 
 \def\sh{\ifmmode\let\next\sh@\else
 \def\next{\errmessage{Use \string\sh\space only in math
 mode}}\fi\next} \def\sh@#1{{\sh@@{#1}}}
 \def\sh@@#1{\fam\eusmfam#1}

\ifcase\@ptsize \font\tenmsa=msam10 \font\sevenmsa=msam7
 \font\fivemsa=msam5 \font\tenmsb=msbm10
 \font\sevenmsb=msbm7 \font\fivemsb=msbm5 \or
 \font\tenmsa=msam10 scaled \magstephalf
 \font\sevenmsa=msam7 \font\fivemsa=msam5
 \font\tenmsb=msbm10 scaled \magstephalf
 \font\sevenmsb=msbm7 \font\fivemsb=msbm5 \or
 \font\tenmsa=msam10 scaled \magstep1 \font\sevenmsa=msam7
 \font\fivemsa=msam5 \font\tenmsb=msbm10 scaled \magstep1
 \font\sevenmsb=msbm7 \font\fivemsb=msbm5 \fi

\newfam\msafam \newfam\msbfam \textfont\msafam=\tenmsa
\scriptfont\msafam=\sevenmsa \scriptscriptfont\msafam=\fivemsa
\textfont\msbfam=\tenmsb \scriptfont\msbfam=\sevenmsb
\scriptscriptfont\msbfam=\fivemsb

\def\Bbb{\ifmmode\let\next\Bbb@\else
 \def\next{\errmessage{Use \string\Bbb\space only in math
 mode}}\fi\next} \def\Bbb@#1{{\Bbb@@{#1}}}
 \def\Bbb@@#1{\fam\msbfam#1} \def\hexnumber@#1{\ifnum#1<10
 \number#1\else \ifnum#1=10 A\else\ifnum#1=11
 B\else\ifnum#1=12 C\else \ifnum#1=13 D\else\ifnum#1=14
 E\else\ifnum#1=15 F\fi\fi\fi\fi\fi\fi\fi}
 \def\msa@{\hexnumber@\msafam} \def\msb@{\hexnumber@\msbfam}
 \mathchardef\square="0\msa@03

\makeatother
\newcommand{\beq}{\begin{equation}}
\newcommand{\eeq}{\end{equation}}
\newcommand{\ba}{\begin{array}}
\newcommand{\ea}{\end{array}}
\newcommand{\bea}{\begin{eqnarray}}
\newcommand{\eea}{\end{eqnarray}}
\newcommand{\bean}{\begin{eqnarray*}}
\newcommand{\eean}{\end{eqnarray*}}
\def\beqa{\begin{eqnarray}}
\def\eeqa{\end{eqnarray}}
\newtheorem{theorem}{Theorem}[section]

\newtheorem{remark}[theorem]{Remark}

\newtheorem{proof}{Proof.}

\makeatother
 \newcommand{\RR}{{\Bbb R}}
\newcommand{\CC}{{\Bbb C}} 
\newcommand{\ZZ}{{\Bbb Z}} \newcommand{\QQ}{{\Bbb Q}}
 \newcommand{\NN}{{\Bbb N}}
\newcommand{\II}{{\Bbb I}}
\begin{document}

\begin{titlepage}

\hfill{DFPD01/TH/18}

\hfill{math.AG/0105051}

\vspace{1cm}
\begin{center}
{\large \bf Eigenfunctions of the Laplacian Acting on Degree Zero
Bundles\\}

\vspace{.3cm}

{\large \bf over Special Riemann Surfaces}

\end{center}
\vspace{1.5cm} \centerline{MARCO MATONE} \vspace{0.8cm}
\centerline{\it Department of Physics ``G. Galilei'' - Istituto
Nazionale di Fisica Nucleare} \centerline{\it University of
Padova} \centerline{\it Via Marzolo, 8 - 35131 Padova, Italy}
\centerline{matone@pd.infn.it}

\vspace{2cm} \centerline{\sc ABSTRACT}

\vspace{0.6cm} \noindent We find an infinite set of eigenfunctions
for the Laplacian with respect to a flat metric with conical
singularities and acting on degree zero bundles over special
Riemann surfaces of genus greater than one. These special surfaces
correspond to Riemann period matrices satisfying a set of
equations which lead to a number theoretical problem. It turns out
that these surfaces precisely correspond to branched covering of
the torus. This reflects in a Jacobian with a particular kind of
complex multiplication.

\end{titlepage}

\newpage

\setcounter{footnote}{0}

\renewcommand{\thefootnote}{\arabic{footnote}}

\renewcommand{\theequation}{\thesection.\arabic{equation}}
\newcommand{\mysection}[1]{\setcounter{equation}{0}\section{#1}}

\mysection{Introduction}

In \cite{M} it was considered the problem of constructing a set of
eigenfunctions of the Laplacian acting on degree zero bundles over
Riemann surfaces. The corresponding metric is given by the modulo
square of a particular holomorphic one-differential
$\omega_{n,m}$. In \cite{M} it has been also shown that
eigenvalues with a nontrivial dependence on the complex structure
may be obtained as solutions of the equation \beq
\omega_{n',m'}=c\,\omega_{n,m}. \label{equazione}\eeq This
equation is equivalent to
\begin{equation}
m_j'-\sum_{k=1}^h\Omega_{jk}n_k'=\overline c
\left(m_j-\sum_{k=1}^h\Omega_{jk}n_k\right),\qquad j=1,\ldots,h,
\label{iudqwwwwINIZIO}\end{equation} where $m_j,n_j,m_j',n_j'$ are
integers. In this paper we find a set of solutions of such an
equation. The general problem involved in Eq.(\ref{equazione})
concerns the properties of the Riemann period matrix and its
number theoretical structure. We will call special Riemann
surfaces those with a Riemann period matrix satisfying
Eq.(\ref{iudqwwwwINIZIO}).

\mysection{Primitive differentials and scalar products}

In this section, after fixing the notation, we introduce an
infinite set of holomorphic one-differentials (we shall call them
{\it primitive differentials}) which can be considered as the
building blocks for our investigation. These differentials have
been previously introduced in \cite{M}. Next, we consider scalar
products defined in terms of monodromy factors associated to the
primitive differentials. By Riemann bilinear relations we
generalize a result in \cite{M} obtaining a relation between
scalar product, monodromy factors and surface integrals. These
aspects are reminiscent of the well-known relations between area,
Fuchsian dilatations and Laplacian spectra which arise for example
in the framework of the Selberg trace formula.

\subsection{Notation and definitions}

Let $\Sigma$ be a compact Riemann surface of genus $h\geq 1$ and
$\alpha_1,\ldots,\alpha_h,\beta_1,\ldots\beta_h$ a symplectic
basis of the first homology group $H_1(\Sigma,{\ZZ})$, that is
with intersection matrix
\begin{equation}
\pmatrix{\alpha\cdot \alpha &\alpha \cdot \beta \cr
\beta\cdot\alpha & \beta\cdot\beta \cr}= \pmatrix{0 & \II \cr -\II
& 0 \cr}, \label{aantrsctnmtrx}\end{equation} where $\II$ is the
$h\times h$ identity matrix. Let $\omega_1,\ldots,\omega_h$ be the
basis of the $\CC$ vector space of holomorphic one-differentials
with the standard normalization
\begin{equation}
\oint_{\alpha_j}\omega_k=\delta_{jk}.
\label{stndnorm}\end{equation} The Riemann period matrix is
defined by
\begin{equation}
\Omega_{jk}\doteq
\oint_{\beta_j}\omega_k.\label{jwelmk}\end{equation}
By means of the Riemann bilinear relations \cite{FK} \beq
\int_{\Sigma}\omega_j\wedge\overline\omega_k=\sum_{l=1}^h\left(
\oint_{\alpha_l}\omega_j\oint_{\beta_l}\overline\omega_k
-\oint_{\alpha_l}\overline\omega_k\oint_{\beta_l}\omega_j\right),
\label{hdhd}\eeq it follows that $\Omega_{ij}=\Omega_{ji}$, and
$\det \Omega^{(2)}>0$, where $\Omega_{kj}^{(1)}\doteq{\Re}\,
\Omega_{kj}$ and $\Omega_{kj}^{(2)}\doteq {\Im}\, \Omega_{kj}$. We
shall denote elements in $H_1(\Sigma,{\ZZ})$ by
\begin{equation}
\gamma_{p,q}\doteq {p}\cdot{\alpha}+{q}\cdot {\beta}, \qquad
(q,p)\in {\ZZ}^{2h},\label{babab}\end{equation} in particular
\begin{equation}
\oint_{\gamma_{p,q}}\omega\doteq\sum_{k=1}^h
\left(p_k\oint_{\alpha_k}\omega+q_k\oint_{\beta_k}\omega \right),
\label{ABZ1}\eeq
 and
\begin{equation}
\int^{z+\gamma_{p,q}}_{z_0}\omega\doteq\int^z_{z_0}\omega
+\oint_{\gamma_{p,q}}\omega,\label{ABZ2}\eeq where,
$z_0\in\Sigma$, $z\in\Sigma$ and $\omega$ is an arbitrary
meromorphic one-differential.

\subsection{Primitive differentials}

Let us consider the set ${\cal K}(\Sigma)\subset H^1(\Sigma)$ of
all nonzero real, harmonic one-forms on $\Sigma$ which are
integral, $i.e.$ such that \beq\oint_\gamma \alpha\in \ZZ,\qquad
\gamma\in H_1(\Sigma,\ZZ).\label{acca1}\eeq Note that ${\cal
K}(\Sigma)$ is a lattice in the $2h$ real vector space
$H^1(\Sigma)$. Let us consider the set of holomorphic
one-differentials \beq {\cal H}(\Sigma)\doteq\{\omega=\pi
i(\alpha+i^\star \alpha)|\alpha\in{\cal
K}(\Sigma)\},\label{iprimitives}\eeq where ${}^\star$ denotes the
conjugation operator whose action on a one-form
$\eta=u(z)dz+v(z)d\bar z$ is ${}^\star\eta=-i u(z)dz+iv(z)d\bar
z$. Note that $\alpha={\Im}\,\omega/\pi$ and ${}^\star
\alpha=-i{\Re}\,\omega/\pi$.

The elements of ${\cal H}(\Sigma)$ correspond to the following
{\it primitive differentials} \cite{M}
\begin{equation}
\omega_{n,m}\doteq\sum_{k=1}^hc_{n,m;k}\omega_k,
\label{n21}\end{equation} where
\begin{equation}
c_{n,m;k}=\pi\sum_{j=1}^h\left(m_j-\sum_{l=1}^hn_l\overline
\Omega_{lj} \right)\left({\Omega^{(2)}}^{-1}\right)_{jk}, \qquad
(n,m)\in {\ZZ}^{2h},\label{gduyg}\end{equation} $k=1,\ldots,h$.
Let us set
\begin{equation}
f_{n,m}(z)\doteq e^{\int^z_{z_0} \omega_{n,m}},
\label{iudqlk}\end{equation} where $z_0$ is fixed on $\Sigma$ and
$z\in\Sigma$. Note that the monodromy of $f_{n,m}$ takes real
values, that is
\begin{equation}
(n,m|q,p)\in{\RR},\qquad (q,p)\in{\ZZ}^{2h}.
\label{iabudqlk}\end{equation} where
\begin{equation}
(n,m|q,p)\doteq e^{\oint_{\gamma_{p,q}}\omega_{n,m}}=
{f_{n,m}(z+\gamma_{p,q})\over
f_{n,m}(z)}=\exp\left[{\sum_{j=1}^h(p_j+\sum_{k=1}^h
q_k\Omega_{kj})c_{n,m;j}}\right]. \label{oihq3}\end{equation} For
later use we define the coefficients
\begin{equation}
D_{kj}^{nm}\doteq m_k\delta_{kj}-n_k\overline\Omega_{kj},
\label{podpl}\end{equation} so that
\begin{equation}
c_{n,m;k}=\pi\sum_{j,l=1}^hD_{jl}^{nm}({\Omega^{(2)}}^{-1})_{lk}.
\label{oidhgf4f}\end{equation}

\subsection{Monodromy and scalar products}

Let us define the scalar product
\begin{equation}
\langle\langle
n,m|q,p\rangle\rangle\doteq\oint_{\gamma_{p,q}}\omega_{n,m},\qquad
(n,m;q,p)\in{\ZZ}^{4h}. \label{iugnx}\end{equation} By
(\ref{gduyg}) we have
\begin{equation}
\langle\langle n,m|q,p\rangle\rangle=\pi\sum_{j,k=1}^h
\left(p_j+\sum_{i=1}^h q_i\Omega_{ij}\right)
\left({\Omega^{(2)}}^{-1}\right)_{jk} \left(m_k-\sum_{l=1}^h
\overline \Omega_{kl}n_l\right), \label{oidjgft}\end{equation}
which has the properties
\begin{equation}
\overline{\langle\langle n,m| q,p\rangle\rangle}= \langle\langle
-q,p|-n,m\rangle\rangle =\langle\langle n,m| q,p\rangle\rangle
-2i\pi(p\cdot n+q \cdot m), \label{ixhuq}\end{equation} so that
\begin{equation}
\Im\,{\langle\langle n,m| q,p\rangle\rangle}= \Im
\,{\langle\langle m,n| p,q\rangle\rangle}.
\label{aixahuq}\end{equation} Furthermore \beq \langle\langle
n,m|q,p\rangle\rangle= {1\over 2\pi i}\sum_{j=1}^h
\left(\langle\langle n,m|\hat\jmath,0\rangle\rangle \langle\langle
0,\hat\jmath| q,p\rangle\rangle +\langle\langle
n,m|0,\hat\jmath\rangle\rangle \langle\langle \hat\jmath,0|
q,p\rangle\rangle\right), \label{gdft}\eeq where $\hat \jmath
\doteq(\delta_{1j},\delta_{2j},\ldots,\delta_{hj})$. By
(\ref{iabudqlk}) and (\ref{ixhuq}) it follows that
\begin{equation}
(n,m|q,p)=e^{\langle\langle n,m| q,p\rangle\rangle}=
e^{\overline{\langle\langle n,m| q,p\rangle\rangle}}= (-q,p|-n,m).
\label{aosdsi}\end{equation} Note that by (\ref{ixhuq}) the scalar
product
\begin{equation}
\langle n,m|q,p\rangle\doteq
{\Re}\,\oint_{\gamma_{p,-q}}\omega_{n,m}= {\Re}\, \langle\langle
n,m| -q,p\rangle\rangle=-{\Re}\, \langle\langle n,m|
q,-p\rangle\rangle, \label{pdfoijw}\end{equation} has the symmetry
property
\begin{equation}
\langle n,m| q,p\rangle= \langle q,p| n,m\rangle,
\label{duiohw}\end{equation} and \beq \langle n,m|
q,p\rangle=\langle\langle n,m|- q,p\rangle\rangle -i\pi (p\cdot n
- q\cdot m). \label{hfgty}\eeq Furthermore, by (\ref{oidjgft}) and
(\ref{ixhuq})
\begin{equation}
\langle n,m|q,p\rangle=\pi\sum_{j,k=1}^h\left[
\left(p_j-\sum_{i=1}^h q_i\Omega_{ij}^{(1)}\right)
\left({\Omega^{(2)}}^{-1}\right)_{jk} \left(m_k-\sum_{l=1}^h
\Omega_{kl}^{(1)}n_l\right) + q_j\Omega^{(2)}_{jk}n_k\right].
\label{oidjgftbb}\end{equation} Positivity of $\Omega^{(2)}$
implies that
\begin{equation}
\langle n,m| n,m\rangle\ge 0, \label{duiohwb}\end{equation} and
$\langle n,m| n,m\rangle=0$ iff $|n,m\rangle =|0,0\rangle$.
Observe that by (\ref{hfgty}) it follows that \beq \langle n,m|
n,m\rangle = \langle\langle n,m| -n,m\rangle\rangle.
\label{gdtey}\eeq $\langle n,m|q,p\rangle$ can be expressed in
terms of the coefficients $c_{n,m;k}$ only \beqa\langle
n,m|q,p\rangle &=&{\pi}^{-1}\sum_{j,k=1}^h
\left(a_{n,m;j}\Omega_{jk}^{(2)}a_{q,p;k}+
b_{q,p;j}\Omega_{jk}^{(2)}b_{n,m;k}\right)\nonumber\\
&=&{\pi}^{-1} \sum_{j,k=1}^h \left( \langle n,m|0,\hat
\jmath\rangle \Omega_{jk}^{(2)} \langle 0,\hat k|q,p\rangle
+q_j\Omega_{jk}^{(2)}n_k \right), \label{interms}\eeqa where
$a_{n,m;k}\doteq {\Re}\,c_{n,m;k}$ and $b_{n,m;k}\doteq
{\Im}\,c_{n,m;k}$, that is
\begin{equation}
a_{n,m;k}=\pi\sum_{j=1}^h\left( m_j-\sum_{l=1}^hn_l
\Omega_{lj}^{(1)}\right)\left({\Omega^{(2)}}^{-1}\right)_{jk},\qquad
b_{n,m;k}=\pi n_k, \label{ABZ5}\eeq $k=1,\ldots,h$, $(n,m)\in
{\ZZ}^{2h}$. Finally, we note that in deriving Eq.(\ref{interms})
we used the relation
$\omega_{n,m}=\sum_{k=1}^h\left(\oint_{\alpha_k}\omega_{n,m}\right)
\omega_k$, that is
\begin{equation}
c_{n,m;k}=\oint_{\alpha_k}\omega_{n,m}=\langle\langle n,m|0,\hat
k\rangle\rangle,\qquad a_{n,m;k}=\langle n,m|0,\hat k\rangle.
\label{iudifdgf}\end{equation}

\subsection{Duality, surface integrals and monodromy}

Let us set \beq \rho_{n,m}^{(1)}\doteq \Re\,\sum_{k=1}^h
d_{n,m;k}^{(1)}\omega_k. \label{gfhty}\eeq We consider the problem
of finding the structure of the coefficients $d_{n,m;k}^{(1)}$
such that \beq \oint_{\gamma_{p,q}}\rho_{n,m}^{(1)}=
\oint_{\gamma_{m,n}}\rho_{q,p}^{(1)}, \label{agdft}\eeq
$(n,m;q,p)\in\ZZ^{4h}$. This equation implies that \beq
d_{n,m;k}^{(1)}=\sum_{j=1}^h
\left(m_j+\sum_{l=1}^hn_l\Omega_{lj}^{(1)}\right)E_{jk}-
\sum_{j=1}^h\left(im_j+\sum_{l=1}^h
n_l\Omega_{lj}^{(2)}\right)F_{jk}+i\sum_{l=1}^hn_lG_{lk},
\label{aaagfthdy}\eeq with $E_{ij}, F_{ij}, G_{ij}$ real symmetric
matrices not depending on $(n,m)$. Note that the condition
$\overline {d^{(1)}_{n,m;k}}= d_{n,m;k}^{(1)}-2\pi i n_k$ is
equivalent to $F_{jk}=0, G_{jk}=\pi\delta_{jk}$. We now set \beq
\rho_{n,m}^{(2)}\doteq \Re\,\sum_{k=1}^h d_{n,m;k}^{(2)}\omega_k,
\label{ribisgfhty}\eeq and consider the condition \beq
\oint_{\gamma_{p,-q}}\rho_{n,m}^{(2)}=
\oint_{\gamma_{m,-n}}\rho_{q,p}^{(2)}, \label{ribisagdft}\eeq that
is \beq d_{n,m;k}^{(2)}=\sum_{j=1}^h
\left(m_j-\sum_{l=1}^hn_l\Omega_{lj}^{(1)}\right)E_{jk}+
\sum_{j=1}^h\left(im_j+\sum_{l=1}^h
n_l\Omega_{lj}^{(2)}\right)F_{jk}+i\sum_{l=1}^hn_lG_{lk}
=\overline{d^{(1)}_{-n,m;k}}. \label{gfthdy}\eeq In this case the
conditions \beq \overline {d_{n,m;k}^{(2)}}= d_{n,m;k}^{(2)}-2\pi
i n_k,\qquad E_{jk}=\pi\left({\Omega^{(2)}}^{-1}\right)_{jk},
\label{hgyflko}\eeq give $d_{n,m;k}^{(2)}=c_{n,m;k}$, that is
\beq\rho_{n,m}^{(2)}=\Re\,\omega_{n,m}. \label{podsfght}\eeq Thus
$c_{n,m;k}$ can be fixed by imposing either the singlevaluedness
of $\exp\Im \,\int^z\omega_{n,m}$ Eq.(\ref{iabudqlk}), or the
duality condition (\ref{ribisagdft}), the same of
Eq.(\ref{duiohw}) satisfied by $\Re\,\omega_{n,m}$, together with
(\ref{hgyflko}).

We will see that positivity and symmetry of $\langle\cdot|
\cdot\rangle$ are at the basis of the fact that the monodromy of
$f_{n,m}$ under a shift of $z$ around $\gamma_{m,-n}$, is
proportional to the area of the metric $|\omega_{n,m}|^2$. This
metric defines a Laplacian of which $(f_{n,m}/\overline
f_{n,m})^k$, $k\in\ZZ$, are eigenfunctions. These aspects are
reminiscent of the well-known relationships between hyperbolic
dilatations and eigenvalues of the Poincar\'e Laplacian.

By the Riemann bilinear relations we have\beqa
\int_\Sigma\omega_{n,m}\wedge \overline\omega_{q,p}
&=&\sum_{j=1}^h\left( \langle\langle n,m| 0, \hat\jmath \rangle
\rangle \overline{\langle\langle q,p | \hat\jmath,0\rangle\rangle}
-\overline{\langle\langle q,p|0,\hat\jmath\rangle\rangle}
\langle\langle n,m| \hat\jmath,0 \rangle \rangle \right)\nonumber\\
&=&\sum_{j=1}^h\left( \langle\langle n,m| 0, \hat\jmath \rangle
\rangle {\langle\langle -\hat \jmath ,0| -q,p\rangle\rangle}
-{\langle\langle n,m|\hat \jmath,0\rangle\rangle} \langle\langle
0,\hat\jmath| -q,p \rangle \rangle\right), \label{hjdfjeu}\eeqa
where (\ref{ixhuq}) has been used. By (\ref{gdft}) and
(\ref{hfgty}) we obtain \beq {i\over
2}\int_\Sigma\omega_{n,m}\wedge \overline\omega_{q,p} =\pi
\langle\langle n,m |-q,p\rangle\rangle= \pi\langle n,m|
q,p\rangle+i\pi^2 (p\cdot n - q\cdot m), \label{hbdheuy}\eeq so
that the monodromy of $f_{n,m}$ corresponds to a surface integral,
that is
 \beq f_{n,m}(z+\gamma_{p,-q})=e^{{1\over 2\pi i}\int_{\Sigma}
\overline \omega_{q,p}\wedge \omega_{n,m} }f_{n,m}(z).
\label{gfhytu}\eeq By (\ref{ixhuq})(\ref{aixahuq}) and
(\ref{hjdfjeu}) the surface integrals $\int_\Sigma
\omega_{n,m}\wedge \overline\omega_{q,p}$ have the properties \beq
\int_\Sigma\omega_{n,m}\wedge \overline\omega_{q,p}= \int_\Sigma
\omega_{q,p}\wedge \overline\omega_{n,m} +4\pi^2(p\cdot n- q\cdot
m), \label{pokju}\eeq and \beq \Im\, {i\over
2}\int_\Sigma\omega_{n,m}\wedge \overline\omega_{q,p}= \Im
\,{i\over 2}\int_\Sigma\omega_{m,n}\wedge \overline\omega_{p,q}.
\label{pkhuye}\eeq The above structure suggests introducing the
holomorphic one-differentials \beq \eta^{(1)}_j(z)= \pi
\sum_{k=1}^h\left({\Omega^{(2)}}^{-1}\right)_{jk}\omega_{k}(z),
\qquad j=1,\ldots,h, \label{uic}\eeq
\begin{equation}
\eta^{(2)}_j(z)=\pi\sum_{l=1}^h\left[i\delta_{jl}-\sum_{k=1}^h
\Omega_{jk}^{(1)}\left({\Omega^{(2)}}^{-1}\right)_{kl}\right]
\omega_{l}(z),\qquad j=1,\ldots,h. \label{uhgdf}\end{equation} By
(\ref{n21}) and (\ref{gduyg}) it follows that
\begin{equation}
\omega_{n,m}=\sum_{k=1}^h(m_k\eta_k^{(1)}+n_k\eta_k^{(2)}),
\label{oidhgfy}\end{equation} moreover
\begin{equation}
{\Im}\,\oint_{\alpha_k}\eta^{(1)}_j=0, \qquad {\Im}\,
\oint_{\beta_k}\eta^{(1)}_j=\pi\delta_{jk}, \label{ABZ7}\eeq
\vspace{0.3cm}
\begin{equation}
{\Im}\,\oint_{\alpha_k}\eta^{(2)}_j=\pi\delta_{jk}, \qquad {\Im}\,
\oint_{\beta_k}\eta^{(2)}_j=0, \label{ABZ8}\eeq $j,k=1,\ldots,h$.
Let us set
\begin{equation}
g_k^{(j)}(z)\doteq\exp \int^z_{z_0} \eta_k^{(j)},\qquad j=1,2,
\label{aoidhgfy}\end{equation} $k=1,\ldots,h$. The expression of
$f_{n,m}$ in terms of ${g_k^{(j)}}$ has the simple form
\begin{equation}
f_{n,m}=\prod_{k=1}^h{g_k^{(1)}}^{m_k}{g_k^{(2)}}^{n_k}.
\label{iuwgfht}\end{equation} Furthermore, since
$\eta_k^{(2)}=-\sum_{j=1}^h\Omega_{kj}\eta_j$,
it follows that
\begin{equation}
f_{n,m}(z)=\prod_{j,k=1}^h
{g_j^{(2)}(z)}^{D_{kj}^{nm}}=e^{\sum_{j,k=1}^hD_{kj}^{nm}\int^z_{z_0}
\eta_j}, \label{iuwgfhtas}\end{equation} where the coefficients
$D_{kj}^{nm}$ have been defined in (\ref{podpl}).

\mysection{Eigenfunctions}

Let $\alpha$ be an element of ${\cal K}(\Sigma)$ and $\omega=\pi i
(\alpha+i^\star \alpha)$ the corresponding holomorphic
differential in ${\cal H}(\Sigma)$. Let $g$ be the metric on
$\Sigma$ given by the line element \beq ds_g=|\omega|.
\label{oiajmc}\eeq In local coordinates, if $\omega=h(z) dz$, then
\beq ds_g^2=|h(z)dz|^2.\label{oiuvG}\eeq Thus $g$ defines a flat
metric on $\Sigma$ with conical singularities at the $2h-2$ zeroes
of $\omega$. For an account on the geodesic dynamics of such
surfaces see for example \cite{Vorobets}. Let $\Delta_g$ be the
corresponding Laplacian acting on degree zero bundles. In local
coordinates \beq \Delta_g=-|h(z)|^{-2}\partial_z\partial_{\bar z}.
\label{oicjP}\eeq

\vspace{.5cm}

\noindent {\bf Theorem 1.} If $k\in\ZZ$ and $F_k$ is the single
valued function \beq F_k(z)=e^{2\pi i
k\int_{z_0}^z\alpha},\label{teorema31}\eeq then
\begin{equation}
\Delta_gF_k=k^2 F_k. \label{eqautov}\end{equation}

\vspace{.5cm}

\noindent {\sc Proof.} Immediate.

\vspace{.5cm}

\noindent Let us define the single valued functions
\begin{equation}
h_{n,m}\doteq{f_{n, m}\over\overline
f_{n,m}}=e^{\int^z\omega_{n,m}-\overline{\int^z\omega_{n,m}}}.
\label{gfhdy}\end{equation} Note that the functions $F_k$
correspond to $h^k_{n,m}$ for some integer vectors $n$ and $m$.
For $k\in {\ZZ\backslash \{0\}}$ we have
\begin{equation}
\int_{\Sigma}\omega_{n,m}\wedge\overline \omega_{n,m}\exp
k\left(\int^z\omega_{n,m} -{\overline
{\int^z\omega_{n,m}}}\right)=0,\label{oiwhd}\end{equation} which
follows from the fact that the integrand is a total derivative.
Since $\overline F_k=F_k^{-1}$, it follows that the $F_k$ satisfy
the orthonormality relation
\begin{equation}
\int_\Sigma d\mu F_k \overline F_j=\delta_{jk},
\label{ijdgfht}\end{equation} where \beq d\mu\doteq
{\omega\wedge\overline\omega\over
\int_\Sigma\omega\wedge\overline\omega}.
\label{definition}\end{equation}


%

\subsection{Multivaluedness, area and eigenvalues}

The area of $\Sigma$ with respect to the metric
$ds^2=|\omega_{n,m}|^2$ is given by \beq A_{n,m}= {i\over 4}
\int_\Sigma \omega_{n,m}\wedge \overline\omega_{n,m}={\pi\over 2}
\langle n,m | n,m\rangle={\pi^2\over2}(m-n\cdot\Omega)\cdot
{\Omega^{(2)}}^{-1}\cdot(m-n\cdot\overline\Omega) .
\label{jhkgi}\eeq The multivaluedness of $f_{n, m}$ is related to
$A_{n,m}$. In particular, winding around the cycle
$\gamma_{n,-m}=-{m} \cdot {\alpha}+ {n} \cdot {\beta}$, we have
\begin{equation}
{\cal P}_{n,-m}f_{n, m}(z)=e^{-{2 A_{n,m}\over \pi}} f_{n,
m}(z)=e^{-\langle n,m|n,m\rangle} f_{n, m}(z),
\label{trmn}\end{equation} where ${\cal P}_{q,p}$ is the winding
operator
\begin{equation}
{\cal
P}_{q,p}g(z)=g(z+\gamma_{p,q}).\label{wndngprtr}\end{equation}
Comparing (\ref{trmn}) with (\ref{eqautov}) we get the following
relationship connecting dilatations and eigenvalues
\begin{equation}
\lambda_k=-{\pi}\log{ {\cal P}_{n,-m}^{k^2}f_{n, m}(z)\over f_{n,
m}(z)}.\label{eigmulti}\end{equation} Thus we can express the
action of the Laplacian on $h_{n,m}^k=(f_{n,m}/\overline
f_{n,m})^k$ in terms of the winding operator acting on $f_{n,m}$.
This relationship between eigenvalues and multivaluedness is
reminiscent of a similar relation arising between geodesic lengths
(Fuchsian dilatations) and eigenvalues of the Poincar\'e Laplacian
(Selberg trace formula). This is not a surprise. Actually, in the
previous sections we have reproduced in the higher genus case some
of the structures arising in the case of the torus. In particular,
we considered the points $c_{n,m;k}$ for which the imaginary part
of $\sum_k c_{n,m;k}\oint_{\gamma_{p,q}}$ takes values in $\pi
\ZZ$. This is reminiscent of the Poisson summation formula where
it (McKean).

\subsection{Genus one}

Let us denote by $\tau$ the Riemann period matrix in the case of
the torus. We set $\tau^{(1)}\doteq{\Re}\,\tau$ and
$\tau^{(2)}\doteq{\Im}\,\tau$. For $h=1$ we have \beq
\omega_{n,m}=c_{n,m}\omega, \label{orifif}\eeq with $\omega\doteq
\omega_1$ the unique holomorphic one-differential on the torus
such that $\oint_\alpha\omega=1$. By (\ref{gduyg}) we have \beq
c_{n,m}=\pi {\left(m-n\overline\tau \right)\over {\tau^{(2)}}},
\qquad (n,m)\in {\ZZ}^{2}. \label{aagduyg}\end{equation} The
functions \beq h_{n,m}=e^{c_{n,m}\int^z\omega -\overline
c_{n,m}\overline{\int^z \omega}},\qquad ({n,m})\in {\ZZ}^2,
\label{ifud}\eeq coincide with the well-known eigenfunctions of
the Laplacian $\Delta=-2\partial_z\partial_{\bar z}$ \beq \Delta
h_{n,m}=\lambda_{n,m} h_{n,m}, \qquad (n,m)\in {\ZZ}^2,
\label{ciugw}\eeq where \beq \lambda_{n,m}\doteq2|c_{n,m}|^2=
{2\pi^2}{(m-n\tau)(m-n\overline\tau)\over{\tau^{(2)}}^2}.
\label{eigtorus}\eeq Under modular transformations of $\tau$ the
eigenvalues transform in the following way \beq
\lambda_{n,m}(\tau+1)=\lambda_{n,m-n}(\tau), \qquad
\lambda_{n,m}\left(-{1\over \tau}\right)=
|\tau|^2\lambda_{-m,n}(\tau). \label{hfgytrr}\eeq The results in
the previous sections shown that the eigenvalues can be generated
by winding around the homology cycles, namely \beq
\lambda_{n,m}=2\pi{\langle n,m|n,m\rangle\over \tau^{(2)}}
 ={2}\pi{{\Re}\,\oint_{\gamma_{m,-n}}
\omega_{n,m}\over \tau^{(2)}}. \label{hdghft}\eeq In other words,
by acting with the winding operator we recover the full spectrum.
We now prove modular invariance of
${\tau^{(2)}}^{-1}{\det'\,\Delta}$ without computing it (the prime
indicates omission of the zero mode of $\Delta$). First of all
note that
\begin{equation}
\mu_{n,m}(\tau)=\tau^{(2)}\lambda_{n,m}(\tau),\quad(n,m)\in \ZZ^2,
\label{ABZ9}\eeq which correspond to the eigenvalues of the
Laplacian $\Delta'=\tau^{(2)}\Delta$, satisfy \beq
\mu_{\gamma(n,m)}(\gamma\cdot\tau)=\mu_{n,m}(\tau),
\label{hfgytrrbis}\eeq where $\gamma(n,m)\doteq(\tilde n,\tilde
m)$ with
\begin{equation}
\left(\begin{array}{c}\tilde m\\ \tilde n
\end{array}\right)
\doteq \left(\begin{array}{c}a\\c
\end{array}\begin{array}{cc}b\\d\end{array}\right)
\left(\begin{array}{c}m\\n
\end{array}\right),
\label{ABZ10}\eeq and
\begin{equation}
\gamma\cdot\tau\doteq{a\tau+b\over c\tau+d},\qquad
\left(\begin{array}{c}a\\c
\end{array}\begin{array}{cc}b\\d\end{array}\right)\in PSL(2,\ZZ).
\label{ABZ11}\eeq By (\ref{hfgytrrbis}) we have
$\mu_{\gamma(n,m)}(\tau)=\mu_{n,m}(\gamma^{-1}\cdot \tau)$, that
is any $\mu_{n,m}$, and therefore $\lambda_{n,m}$, can be obtained
{}from a given eigenvalue by modular transformations. The
determinants of $\Delta$ and $\Delta'$ are related by
\begin{equation}
{\det}'\,\Delta'={\det}'\,( \tau^{(2)}){\det}'\,\Delta=
{\tau^{(2)}}^{-1}{\det}\,( \tau^{(2)}){\det}'\,\Delta,
\label{ABZ12}\eeq where we used the fact that $c \, {\det}'\,
c={\det}\, c$, $c\in \CC$. Since ${\det}\, \tau^{(2)}$ can be
regularized by standard techniques, e.g. by the $\zeta$-function
regularization method, to a finite $\tau$-independent constant, we
have
\begin{equation}
{{\det}'\,\Delta\over {\tau^{(2)}}} =const
\prod_{(n,m)\in\ZZ^2\backslash\{0,0\}}\mu_{n,m} =const
\prod_{(n,m)\in\ZZ^2\backslash\{0,0\}}\tau^{(2)}\lambda_{n,m}.
\label{ABZ13}\eeq Modular invariance of
${\tau^{(2)}}^{-1}{\det}'\,\Delta$ now follows by
(\ref{hfgytrrbis}), namely \beq
\prod_{(n,m)\in\ZZ^2\backslash\{0,0\}}\tau^{(2)}\lambda_{n,m}(\tau)=
\prod_{\gamma\in PSL(2,\ZZ)}\mu_{N,M}(\gamma\cdot\tau),
\label{hahahaha}\eeq where $N,M$ are arbitrary integers not
simultaneously vanishing. It follows that ${\det}'\,\Delta$ can be
expressed as the product of all modular transformations acting on
an arbitrary eigenvalue \beq {\det}'\,\Delta=const\, {\tau^{(2)}}
\prod_{\gamma\in PSL(2,\ZZ)}\mu_{N,M}(\gamma\cdot\tau).
\label{ashahahaha}\eeq Note that modular invariance of
${\tau^{(2)}}^{-1}{\det'\,\Delta}$ essentially implies that
${\tau^{(2)}}^{-1}{\det'\,\Delta}=\tau^{(2)}|\eta(\tau)|^4$, where
$\eta(\tau)$ is the Dedekind $\eta$-function.

\mysection{Special Riemann surfaces}

In general there are other eigenfunctions besides $h_{n,m}^k$,
$k\in {\NN}$. For example, when all the $2m_j$'s and $2n_j$'s are
integer multiple of an integer $N$ then the eigenfunctions include
$h_{n,m}^k$, $k\in {\NN}$ whose eigenvalue is $2A_{n,m}k^2/N^2$.

More generally one should investigate whether the period matrix
has some non trivial number theoretical structure. To see this we
first note the trivial fact that since for $h=1$ the space of
holomorphic one-differentials is one-dimensional, it follows that
the ratio between $\omega_{n,m}$ and $\omega_{n',m'}$ is always a
constant. This allows one to construct the infinite set of
eigenvalues labelled by two integers $(n,m)\in{\ZZ}^2$. In the
case $h\ge 2$ the ratio $\omega_{n,m}/\omega_{n',m'}$ is in
general not a constant. This is the reason why we considered the
eigenfunctions of the kind $h_{n,m}^k$ with fixed $(n,m)\in
{\ZZ}^{2h}$. However there are other interesting possibilities.
For example, if besides $\alpha$ also ${}^\star\alpha$ is
integral, then for $k_1,k_2\in\ZZ$, the single valued function
\beq F_{k_1,k_2}= e^{2\pi i\int_{z_0}^z
(k_1\alpha+k_2{}^\star\alpha)},\label{esempio}\eeq satisfies \beq
\Delta_g F_{k_1,k_2} =(k_1^2+k_2^2)F_{k_1,k_2}.
\label{esempios}\eeq Note that since
$k_1\alpha+k_2{}^\star\alpha=2i{\Im}\,[(k_1-ik_2)\omega]$, it
follows that integrality of both $\alpha$ and ${}^\star\alpha$ is
a particular case of a more general one.

\vspace{.5cm}

\noindent {\bf Theorem 2.} If the holomorphic one-differentials
$\omega_{n,m}$ and $\omega_{n',m'}$ satisfy the equation
\begin{equation}
\omega_{n',m'}(z)=c\,\omega_{n,m}(z),
\label{newwsad}\end{equation} with both $(n,m)$ and $(n',m')$ in
${\ZZ}^{2h}$ and $c\in \CC\backslash \QQ$, then the single valued
function
\begin{equation}
h_{n', m'}=e^{c\int^z\omega_{n,m}-\overline c \overline
{\int^z\omega_{n,m}}}\neq {h_{n, m}^k}\qquad c\in {\CC}\backslash
{\QQ},\quad k\in {\QQ}, \label{a1a2e3i4g5s}\end{equation}
satisfies
\begin{equation}
\Delta_{g^{n,m}}h_{n',m'}=\lambda_c h_{n',m'},
\label{a1a2e3i4g5s6787}\end{equation} where \beq
\lambda_c=2A_{n,m}|c|^2, \label{eigencc}\eeq and
\begin{equation}
c={m_i'-\sum_{k=1}^h\overline \Omega_{ik}n_k' \over
m_i-\sum_{k=1}^h\overline \Omega_{ik}n_k}
={m_j'-\sum_{k=1}^h\overline \Omega_{jk}n_k' \over
m_j-\sum_{k=1}^h\overline \Omega_{jk}n_k},\quad i,j=1,\ldots,h.
\label{iudqdpolk}\end{equation}

\vspace{.5cm}

\noindent {\sc Proof.} It is immediate to check
Eq.(\ref{a1a2e3i4g5s6787}). The only point is to find the
expression of $c$. This follows by the observation that since the
holomorphic one-differentials $\omega_1,\ldots , \omega_h$ are
linearly independent, Eq.(\ref{newwsad}) is equivalent to
\begin{equation}
m_j'-\sum_{k=1}^h\Omega_{jk}n_k'=\overline c
\left(m_j-\sum_{k=1}^h\Omega_{jk}n_k\right),\qquad j=1,\ldots,h.
\label{iudqwwww}\end{equation}

\vspace{.5cm}


Note that to each $(n,m)$ and $(n',m')$ satisfying
(\ref{iudqwwww}) there is a possible value of $c\equiv
c(n,m;n',m')$.

Thus we can reproduce in higher genus the basic structure
Eqs.(\ref{ifud})-(\ref{eigtorus}) considered in the torus case and
then to find eigenvalues with a non trivial dependence on the
complex structure.

The problem of finding the possible (in general complex) solutions
of (\ref{newwsad}) is strictly related to the number theoretical
properties of $\Omega$.




Observe that \beq \Delta_{g^{n',m'}}h_{n,m}={\lambda_c}' h_{n,m},
\label{a1a2e3i4g5s6787i}\end{equation} where \beq
{\lambda_c}'=4{A_{n,m}A_{n',m'}\over \lambda_c}.
\label{gfhtyru}\eeq We will call {\it Special}, the Riemann
surfaces admitting solutions of (\ref{iudqwwww}) with non rational
values of $c$. Since a change in the sign of $c$ is equivalent to
a change sign of either $(m,n)$ or $(m',n')$, in the following,
without loss of generality, we will assume that \beq \Im\, \bar
c>0.\label{iasuhbx}\eeq

\subsection{The solution space}

Eq.(\ref{iudqwwww}) has solutions if there are integers $(n,m)$
and $(n',m')$ such that the period matrix satisfies the $h-1$
consistency conditions
\begin{equation}
{m_i'-\sum_{k=1}^h\Omega_{ik}n_k' \over
m_i-\sum_{k=1}^h\Omega_{ik}n_k}= {m_j'-\sum_{k=1}^h\Omega_{jk}n_k'
\over m_j-\sum_{k=1}^h\Omega_{jk}n_k},\quad i,j=1,\ldots,h.
\label{aiaudqdpolk}\end{equation} Thus, for $(n,m)\in \ZZ^{2h}$
fixed, the Laplacian $\Delta_{g^{n,m}}$ has eigenvalues
\begin{equation}
\lambda_c=2A_{n,m}\left|{m_i'-\sum_{k=1}^h\overline
\Omega_{ik}n_k' \over m_i-\sum_{k=1}^h\overline
\Omega_{ik}n_k}\right|^2, \qquad (n',m')\in {\cal
S}_{n,m}(\Omega), \label{aiauadqdpolk}\end{equation} where ${\cal
S}_{n,m}(\Omega)$ denotes the {\it solution space} \beq {\cal
S}_{n,m}(\Omega)\doteq\left\{(n',m')\left|\omega_{n',m'}=c\,
\omega_{n,m}, (n',m')\in \ZZ^{2h}\right\} \right. .
\label{kjhuyi}\eeq For a given $\Omega$ the space ${\cal
S}_{n,m}(\Omega)$ may contain other points besides $(kn,km),
k\in\ZZ$. To investigate the structure of such space we should
understand the nature of the Riemann surfaces whose period matrix
satisfies Eq.(\ref{iudqwwww}). To this end we note that
Eq.(\ref{iudqwwww}) has been suggested in order to reproduce in
higher genus the basic structure Eqs.(\ref{ifud})-(\ref{eigtorus})
considered in the torus case. So, we should expect a Riemann
surface strictly related to the torus geometry. Remarkably, this
is in fact the case as we have the following

\vspace{.5cm}

\noindent {\bf Theorem 3.} The Riemann surfaces with period
matrices satisfying Eq.(\ref{iudqwwww}) correspond to branched
covering of the torus.

\vspace{.5cm}

\noindent {\sc Proof.}  First of all note that by (\ref{iudqwwww})
it follows that the function
\begin{equation}
w(z)=\int^z\hat\omega, \label{paraponziponziperoxxx}\end{equation}
where $\hat\omega \doteq(\bar c n-n')\cdot\omega$, has monodromy
\begin{equation}
\oint_{\gamma_{p,q}}\hat\omega=-p\cdot n'-q\cdot m'+ \bar c(p\cdot
n+q\cdot m), \label{paraponziponzipero}\end{equation} implying
that $w$ is a holomorphic map from $\Sigma$ to the torus with
period matrix $\bar c$. Viceversa, if $w$ is a holomorphic map of
a branched covering of the torus to the torus itself, then
$w(z+\gamma_{p,q})=w(z)+p\cdot N'+q\cdot M' +\tau(p\cdot N +
q\cdot M)$, for some integer vectors $M,N,M',N'$, and by the
Riemann bilinear relations \beq 0=\int_\Sigma\omega_k\wedge
dw=M_k'+\tau M_k -\sum_{j=1}^h(N'_j+\tau N_j)\Omega_{jk},\qquad
k=1,\ldots,h, \label{whichisprecisely}\eeq which is
Eq.(\ref{iudqwwww}) with $\tau=\bar c$ (see \cite{RB} for explicit
constructions of branched covering of the torus\footnote{I am
grateful to the anonymous referee for suggesting Ref.\cite{RB}.}).

\vspace{.5cm}

\noindent {\bf Remark 1.} Eq.(\ref{iudqwwwwINIZIO}), derived in \cite{M}
(see Eqs.(5.18) and (5.23) there), has been subsequently and independently derived in
\cite{G} by studying the null compactification of type-IIA-string perturbation
theory at finite temperature. In \cite{G} an equivalent proof of {\bf Theorem 3.} is also
provided.

\vspace{.5cm}

\subsection{Metric induced by the covering}

By means of the map $w$ from $\Sigma$ to the torus with period
matrix $sign\,({\Im}\,(\bar c))\bar c$ we can construct an
infinite set of eigenfunctions for the Laplacian
$\Delta=-2|\hat\omega|^{-2}\partial_z\partial_{\bar z}$ on
$\Sigma$ defined with respect to the metric $|\hat\omega|^2$.
These eigenfunctions are \beq h_{n,m}=e^{c_{n,m}\int^z\hat\omega
-\overline c_{n,m}\overline{\int^z \hat\omega}},\qquad ({n,m})\in
{\ZZ}^2, \label{ifudbbiiss}\eeq corresponding to the eigenvalues
\beq \lambda_{n,m}=2|c_{n,m}|^2= {2\pi^2}{({m}-c
{n})({m}-\overline c {n})\over {c^{(2)}}^2},\label{eigtorusss}\eeq
where \beq c_{n,m}=\pi {\left(m-n\overline c\right)\over
{c^{(2)}}}, \qquad (n,m)\in
{\ZZ}^{2}.\label{aagduygbbiiss}\end{equation} The proof is a
direct consequence of the relation
$\Delta=-2|\hat\omega|^{-2}\partial_z\partial_{\bar z}=-2
\partial_w\partial_{\bar w}$.

Let us set \beq D_j\doteq\sum_{k=1}^h \overline D_{kj}^{mn},\qquad
D_j'\doteq\sum_{k=1}^h \overline D_{kj}^{m'n'}.\label{oivjf}\eeq
By Eq.(\ref{iudqdpolk}) and the Poisson summation formula, we have
\beq \sum_{n=-\infty}^{\infty}e^{-n^2\pi {D'_j\over D_j}}=
\sqrt{D_j\over D_j'}\sum_{n=-\infty}^{\infty}e^{-n^2\pi {D_j\over
D_j'}},\qquad j=1,\ldots,h.\label{PSF}\eeq

\subsection{Genus 2}

Before considering Eq.(\ref{aiaudqdpolk}) for arbitrary genus, it
is instructive to investigate the $h=2$ case. Since
$\Omega_{ij}=\Omega_{ji}$ we have
\begin{equation}
{m_1'-\Omega_{11}n_1'-\Omega_{12}n_2' \over
m_1-\Omega_{11}n_1-\Omega_{12}n_2}=
{m_2'-\Omega_{12}n_1'-\Omega_{22}n_2' \over
m_2-\Omega_{12}n_1-\Omega_{22}n_2}.
\label{abiabudqdpolk}\end{equation} Thus the problem is the
following: given $\Omega_{11}$ and $\Omega_{12}$ find all the
integers $(n,m;n',m')\in{\ZZ}^{8}$ such that (\ref{abiabudqdpolk})
is satisfied. The solution of this problem depends on structure of
$\Omega$.

We consider period matrices satisfying the relation
\begin{equation}
\Omega_{22}={\hat N_1\over \hat N_4}\Omega_{11}+ {\hat N_2\over
\hat N_4}\Omega_{12}+{\hat N_3\over \hat N_4},\qquad (\hat
N_1,\hat N_2,\hat N_3,\hat N_4)\in {\ZZ}^4,\qquad \hat N_4\ne 0,
\label{m1}\end{equation} that is
\begin{equation}
\Omega_{ij}= \pmatrix{\Omega_{11}&\Omega_{12}\cr
 \Omega_{12}& N_1\Omega_{11}+N_2\Omega_{12}+N_3\cr},
\label{ntrsctnmtrx}\end{equation} where $N_i\doteq \hat N_i/\hat
N_4$, $i=1,2,3$. Positivity of $\Omega^{(2)}_{ij}$ implies the
following condition on $N_1,N_2,N_3$
\begin{equation}
\Omega_{11}^{(2)}\left(N_1 \Omega_{11}^{(2)} + N_2
\Omega_{12}^{(2)}+N_3\right)> {\Omega_{12}^{(2)}}^2.
\label{ifwmdj}\end{equation} With the position (\ref{m1})
Eq.(\ref{abiabudqdpolk}) becomes
\begin{equation}
{m_1'-\Omega_{11}n_1'-\Omega_{12}n_2' \over
m_1-\Omega_{11}n_1-\Omega_{12}n_2}=
{m_2'-N_3n_2'-\Omega_{11}N_1n_2'-\Omega_{12}(n_1'+N_2n_2') \over
m_2-N_3n_2-\Omega_{11}N_1n_2-\Omega_{12}(n_1+N_2n_2)}.
\label{m2}\end{equation} We will look for solutions of this
equation of the form \beq
m_2-N_3n_2-\Omega_{11}N_1n_2-\Omega_{12}(n_1+N_2n_2)
=N\left[m_1-\Omega_{11}n_1-\Omega_{12}n_2\right],
\label{jfhfy}\eeq \beq
m_2'-N_3n_2'-\Omega_{11}N_1n_2'-\Omega_{12}(n_1'+N_2n_2')=
N\left[m_1'-\Omega_{11}n_1'-\Omega_{12}n_2'\right],
\label{lfkio}\eeq with $N\in\QQ$. Since $(n,m;n',m')\in
{\ZZ}^{8}$, the general solutions not depending on $\Omega_{11}$
and $\Omega_{12}$ are
\begin{equation}
m_2-N_3n_2=Nm_1,\quad n_1+N_2n_2=Nn_2, \quad N_1n_2=Nn_1,
\label{m3}\end{equation}
\begin{equation}
m_2'-N_3n_2'=Nm_1',\quad n_1'+N_2n_2'=Nn_2', \quad N_1n_2'=Nn_1'.
\label{m4}\end{equation} Note that each solution of Eq.(\ref{m3})
defines a metric $g^{n,m}$ whereas the solutions of Eq.(\ref{m4})
give, by (\ref{aiauadqdpolk}), the eigenvalues $\lambda_c$ of
$\Delta_{g^{n,m}}$.

The compatibility condition for (\ref{m4}) constrains $N$ to be
\begin{equation}
N_\pm={N_2\pm \sqrt{N_2^2+4N_1}\over 2}. \label{m5}\end{equation}
Since $N\in\QQ$ and $N_i\in {\ZZ}/\hat N_4$, $i=1,2,3$, we have
\begin{equation}
N_1=MN_2+M^2, \qquad M\in \left\{k\in \QQ \left|(kN_2+k^2)\in
\ZZ/\hat N_4 \right\}.\right. \label{m6}\end{equation}
Eq.(\ref{m3}) has a double set of solutions. In the case
$N=N_+=N_2+M$, we have \beq n_2={n_1\over M},\qquad
m_2=(N_2+M)m_1+{N_3n_1\over M},\qquad (n_1,m_1)\in\Gamma^{(+)},
\label{firstset}\eeq where
\begin{equation}
\Gamma^{(+)}\doteq\left\{(k,j)\in\ZZ^2\left|\left({k\over M},
(N_2+M)j+{N_3k\over M}\right)\in \ZZ^2\right\}.\right.
\label{ABZ14}\eeq In the second case $N=N_-=-M$, so that \beq
n_2=-{n_1\over N_2+M},\qquad m_2=-Mm_1-{N_3n_1\over N_2+M},\qquad
(n_1,m_1)\in\Gamma^{(-)}, \label{secondset}\eeq where
\begin{equation}
\Gamma^{(-)}\doteq\left\{(k,j)\in\ZZ^2\left|\left( {k\over N_2+M},
Mj+{N_3k\over N_2+M}\right)\in \ZZ^2\right\}.\right.
\label{ABZ15}\eeq Note that given $N_1$, $N_2$ and $N_3$ we found
that $N_1=MN_2+ M^2$ and either $N=N_2+M$ or $-M$. Therefore it is
natural to choose $M$, $N_2$ and $N_3$ to parametrize
$\Omega_{22}$. Note also that by (\ref{m1}) and (\ref{m6}) it
follows that $M$ and $-N_2-M$ correspond to the same value of
$\Omega_{22}$.

Given $M, N_2, N_3, N_4$ there are two sets of Laplacians
parametrized by points in $\Gamma^{(\pm)}$. Let $(n_1,m_1)$ be a
point in $\Gamma^{(+)}$ and $(n_2,m_2)$ given by (\ref{firstset}).
The first set of eigenvalues is \beq \lambda_c^{(+)}=2A_{n,m}{
\left|Mm_1'-n_1'\left(M\Omega_{11}+{\Omega_{12}}\right)\right|^2
\over \left|Mm_1-n_1\left(M\Omega_{11}
+{\Omega_{12}}\right)\right|^2}, \label{jghytiu}\eeq where
\begin{equation}
n_2'={n_1'\over M},\quad m_2'=(N_2+M)m_1'+{N_3n_1'\over M}, \qquad
(n_1',m_1')\in \Gamma^{(+)}. \label{ABZ16}\eeq Let $(n_1,m_1)$ be
a point in $\Gamma^{(-)}$ and $(n_2,m_2)$ given by
(\ref{secondset}). The second set of eigenvalues is \beq
\lambda_c^{(-)}=2A_{n,m}{ \left|\left(N_2+M\right)m_1'-
n_1'\left[\left(N_2+M\right)\Omega_{11}-{\Omega_{12}}\right]\right|^2
\over \left|(N_2+M)m_1-n_1 \left[\left(N_2+M\right)\Omega_{11}
-{\Omega_{12}}\right]\right|^2}, \label{jghytiasu}\eeq where
\begin{equation}
n_2'=-{n_1'\over N_2+M},\quad m_2'=-Mm_1-{N_3n_1\over N_2+M},
\qquad (n_1',m_1')\in \Gamma^{(-)}. \label{ABZ1ancora}\eeq These
eigenvalues have a structure which is similar to the structure of
the ones of the Laplacian on the torus.

\subsection{Higher genus}

We now consider the higher genus case. By (\ref{aiaudqdpolk}) we
have \beq m_i-\sum_{k=1}^h\Omega_{ik}n_k=N_{ij}\left(
m_j-\sum_{k=1}^h\Omega_{jk}n_k\right), \label{m17}\eeq \beq
m_i'-\sum_{k=1}^h\Omega_{ik}n_k'=N_{ij}\left(
m_j'-\sum_{k=1}^h\Omega_{jk}n_k'\right). \label{m18}\eeq
$i,j=1,\ldots,h$. Eq.(\ref{m17}) gives the constraint on the
structure of the metric $g^{n,m}$. Note that\beq
N_{ij}N_{jk}=N_{ik}, \quad i,j=1,\ldots,h, \label{m20}\eeq in
particular, $N_{ij}=N_{ji}^{-1}$. The matrix $N_{ij}$ is
determined by $h-1$ elements. For example, since
$N_{ij}=N_{i1}N_{1j}=N_{1i}^{-1}N_{1j}$, in terms of
$N_{12},\ldots,N_{1h}$ we have
\begin{equation}
N_{ij}=\pmatrix{1& N_{12} & N_{13} & \ldots & N_{1h}\cr
N_{12}^{-1} & 1 & N_{13}N_{12}^{-1}& \ldots & N_{1h}N_{12}^{-1}\cr
\vdots & \vdots & \vdots & \vdots & \vdots \cr N_{1h}^{-1} &
N_{1h}^{-1}N_{12} & N_{1h}^{-1}N_{13} &\ldots & 1\cr},
\label{m22}\end{equation} which has vanishing determinant.

Since $\Omega$ is symmetric, it follows that in each one of the
$h-1$ equations Eq.(\ref{m17}) (or Eq.(\ref{m18})) there is always
one, and only one, matrix element appearing in both sides. In
other words, both sides of Eq.(\ref{m17}) contain
$\Omega_{ij}=\Omega_{ji}$. Therefore, it is natural to consider
period matrices of the form \beq
\Omega_{ij}=\sum_{k,l=1}^hN_{ij}^{kl}\Omega_{kl}+M_{ij}, \qquad
N_{ij}^{kl}\in\QQ, \quad M_{ij}\in\QQ, \qquad N_{ij}^{ij}=0,
\label{uigwlk}\eeq and then, for each pair $i,j$, to substitute it
in the equation involving $N_{ij}$. This allows us to transform,
for each $i,j$, Eqs.(\ref{m17})(\ref{m18}) in equations containing
all the matrices elements but $\Omega_{ij}$ on both sides. Note
that the symmetry of $\Omega_{ij}$ implies that $N_{ij}^{kl}$ is
symmetric in $i,j$ and $k,l$ separately and $M_{ij}=M_{ji}$. By
(\ref{uigwlk}) we have that Eqs.(\ref{m17})(\ref{m18}) become \beq
N_{ij}m_j-m_i+\sum_{k,l=1}^h\left\{(n_j-n_iN_{ij})
\left(N_{ij}^{kl}\Omega_{kl}+M_{ij}\right)+
\Omega_{kl}\left[\delta_{ik}(n_l-
\delta_{lj}n_j)-\delta_{jk}(n_l-\delta_{li}n_i)N_{ij}\right]\right\}=0,
\label{vdsgd}\eeq \beq
N_{ij}m_j'-m_i'+\sum_{k,l=1}^h\left\{(n_j'-n_i'N_{ij})
\left(N_{ij}^{kl}\Omega_{kl}+M_{ij}\right)+
\Omega_{kl}\left[\delta_{ik}(n_l'-
\delta_{lj}n_j')-\delta_{jk}(n_l'-\delta_{li}n_i')N_{ij}\right]\right\}=0,
\label{vdsgdueee}\eeq $i,j=1,\ldots,h$. We do not investigate the
conditions following from Eqs.(\ref{m17})(\ref{m18}) further,
rather we shortly consider the period matrices satisfying the
conditions \beq
\Omega_{ik}=\sum_{l=1}^hN_{ik,j}^l\Omega_{jl}+M_{ik}, \qquad
N_{ik,j}^l\in \QQ,\; M_{ik}\in \QQ, \quad
i,j,k=1,\ldots,h.\label{m24}\eeq Substituting $\Omega_{ik}$ in the
left hand side of (\ref{m17})(\ref{m18}) these transform in
simplified equations as now they involve matrix elements with the
same value of the first index.

Substituting $\Omega_{jl}$ in the RHS of (\ref{m24}) with
$\sum_{m=1}^hN_{jl,n}^m\Omega_{nm}+M_{jl}$, we have \beq
\Omega_{ik}= \sum_{l=1}^hN_{ik,j}^l\left(\sum_{m=1}^h
N_{jl,n}^m\Omega_{nm}+M_{jl}\right)+M_{ik}, \quad
i,j,k,n=1,\ldots,h. \label{m25}\eeq Comparing (\ref{m24}) with
(\ref{m25}) one obtains a set of equations that, once one makes
the additional requirement that the terms involving the period
matrix cancel separately, become \beq
\sum_{l=1}^hN_{ik,j}^lN_{jl,n}^m=N_{ik,n}^m,\quad
i,j,k,m,n=1,\ldots,h, \label{m27}\eeq \beq \sum_{l=1}^h N_{ik,j}^l
M_{jl}=0,\quad i,j,k=1,\ldots,h. \label{m28}\eeq

\subsection{Special Riemann surfaces and complex multiplication}

We now show that Special Riemann surfaces have a Jacobian with
Complex Multiplication (CM).\footnote{Jacobian with CM have been
recently considered in the framework of rational CFT \cite{GV}.}
First note that in terms of \beq v_k\doteq
m_k-\sum_{j=1}^h\Omega_{kj}n_j,\qquad v_k'\doteq
m_k'-\sum_{j=1}^h\Omega_{kj}n_j',\label{CM1}\eeq
Eq.(\ref{aiaudqdpolk}) reads \beq v_jv_k'-v_kv_j'=0, \qquad
\forall j,k,\eeq which is equivalent to \beq \Omega N\Omega
+\Omega M -\tilde M \Omega -M'=0, \label{CM2}\eeq where $\tilde{}$
denotes the transpose and \beq N_{jk}\doteq n'_jn_k-n_k'n_j,\qquad
M_{jk}\doteq m'_jn_k-m_kn'j,\qquad M'_{jk}\doteq
m'_jm_k-m_k'm_j.\label{CM3}\eeq On the other hand, a Jacobian is
said to admit complex multiplication if there exist integer
matrices $N$, $M$, $M'$ and $N'$ such that \beq \Omega
(M+N\Omega)=M'+N'\Omega,\label{CM4}\eeq that is \beq \Omega
N\Omega+\Omega M-N'\Omega-M'=0.\label{CM5}\eeq Comparing with
(\ref{CM2}), we see that the Jacobians of Special Riemann surfaces
admit a particular kind of CM. According to Theorem 3. this CM is
the one of Jacobian corresponding to branched covering of the
torus.

\end{document}